\documentstyle[12pt,epsf]{article}
\newcommand{\rr}{{\rm{I \! R}}}

\newcommand{\al}{\alpha}
\newcommand{\ds}{\displaystyle}
\newcommand{\f}{\frac}
\newcommand{\lam}{\lambda}

\newcommand{\om}{\omega}
\newcommand{\de}{\delta}

\newcommand{\ca}{{\cal A}}
\newcommand{\th}{\theta}

\title{\bf Hopf bifurcation in a dynamic IS-LM model with time delay}
\author{Mihaela Neam\c tu$^{a}$\thanks{Corresponding author}, Dumitru Opri\c s$^{b}$, Constantin Chil\v arescu$^a$}
\date{}

\begin{document}
\maketitle

\begin{tabular}{cccccccc}
\scriptsize{$^{a}$Department of Economic Informatics and Statistics,}\\
\scriptsize{Faculty of Economics, West University of Timi\c soara,}\\
\scriptsize{str. Pestalozzi, nr. 16A, 300115, Timi\c soara, Romania,}\\
\scriptsize{E-mails:mihaela.neamtu@fse.uvt.ro, cchilarescu@rectorat.uvt.ro}\\
\scriptsize{$^{b}$ Department of Applied
Mathematics, Faculty of Mathematics,}\\
\scriptsize{West University of Timi\c soara, Bd. V. Parvan, nr. 4, 300223, Timi\c soara, Romania,}\\
 \scriptsize{e-mail: opris@math.uvt.ro}\\

\end{tabular}

\begin{quote} \small{{\bf Abstract.} In this paper we investigate the impact of
delayed tax revenues on the fiscal policy out-comes.
Choosing the delay as a bifurcation parameter we study the
direction and the stability of the bifurcating periodic solutions.
With respect to the delay we show when the system is stable. Some
numerical examples are finally given for justifying the
theoretical results.}
\end{quote}

\noindent{\small{{\it Keywords:} delay differential equation,
stability, Hopf bifurcation, IS-LM model.

}}

\section*{\normalsize\bf 1. Introduction}

\hspace{0.6cm}

The differential equations with time delay play an important role
for economy, engineering, biology and social sciences, because a
lot of phenomena are described with their help. In this paper, we
consider a model from economy, of the IS-LM type with time delay
and we study how the delay affects the macroeconomic stability.
The Hopf bifurcation and normal form theories are tools for
establishing the existence and the stability of the periodic
solutions. Similar ideas can be found in [5], [9], [10].

In [2] Cesare and Sportelli taking into account the papers [6],
[7], [8], study a dynamic IS-LM model of the following type:
$$\begin{array}{l}
\dot Y(t)\!=\!\al[I(Y(t),r(t))\!+\!g\!-\!S(Y(t)\!-\!T(Y(t), Y(t\!-\!\tau))))\!-\!T(Y(t), Y(t\!-\!\tau)))]\\
\dot r(t)=\beta[L(Y(t),r(t))-M(t)]\\
\dot M(t)=g-T(Y(t), Y(t\!-\!\tau))),\end{array}$$ with Y as
income, I as investment, g as government expenditure (constant), S
as savings, T as tax revenues, r as rate of interest, L as
liquidity, M as real money supply and $\al, \beta$ as positive
constants. The time delay $\tau$ appears in function T:
$$T(Y(t), Y(t\!-\!\tau)))=d(1-\varepsilon)Y(t)+d\varepsilon Y(t-\tau),\eqno(1)$$
where $d\in (0,1)$ is a common average tax rate and
$\varepsilon\in (0, 1)$ is the income tax share.

Based on the papers [2], [1], [11], we consider the following
IS-LM model:
$$\begin{array}{l}
\dot Y(t)\!\!=\!\!\al[I(Y(t),r(t))\!+\!g\!-\!S(Y(t)\!\!-\!\!T(Y(t), Y(t\!-\!\tau)))\!-\!T(Y(t), Y(t\!-\!\tau))]\\
\dot r(t)=\beta[L(Y(t), r(t))-M(t)]\\
\dot K(t)=I(Y(t-\tau), r(t))-\de K(t)\\
\dot M(t)=g-T(Y(t), Y(t\!-\!\tau)),\end{array}\eqno(2)$$ with
$\delta>0$ and with the initial conditions: $$Y(t)=\varphi(t),
t\in[-\tau, 0], r(0)=r_1, M(0)=M_1, K(0)=K_1.$$

In the following analysis we will consider the function T given by
(1), the investment, the saving  and the liquidity of the form:
$$\begin{array}{l}I(Y(t), r(t))=aY(t)^{\al_1}r(t)^{-\al_2},\quad a>0, \al_1>0,
\al_2>0,\\
S(Y(t)-T(Y(t), Y(t-\tau)))=s(Y(t)-T(Y(t), Y(t-\tau))),\quad s\in(0,1),\\
L(r(t))=mY(t)+\ds\f{\gamma_0 }{r(t)-r_2}, \quad m>0, \gamma_0>0,
r_2>0.\end{array}\eqno(3)$$

The paper is organized as follows. In section 2 we investigate the
local stability of the equilibrium point associated to system (2).
Choosing the delay as a bifurcation parameter some sufficient
conditions for the existence of Hopf bifurcation are found. In
section 3 there is the main aim of the paper, namely the
direction, the stability and the period of a limit cycle solution.
Section 4 gives some numerical simulations which show the
existence and the nature of the periodic solutions. Finally, some
conclusions are given.

\section*{\normalsize\bf 2. The qualitative analysis of system (2).}

\hspace{0.6cm}

Using functions (1) and (3), system (2) becomes:
$$\begin{array}{l}
\dot Y(t)\!=\!\al[((s\!-\!1)(1\!-\!\varepsilon)d\!-\!s)Y(t)\!+\!d\varepsilon(s\!-\!1)Y(t\!-\!\tau)\!+\!aY(t)^{\al_1}r(t)^{-\al_2}\!+\!g]\\
\dot r(t)=\beta[mY(t)+\ds\f{\gamma_0}{r(t)-r_2}-M(t)]\\
\dot K(t)=aY(t-\tau)^{\al_1}r(t)^{-\al_2}-\de K(t)\\
\dot M(t)=g-d(1-\varepsilon)Y(t)-d\varepsilon
Y(t-\tau),\end{array}\eqno(4)$$ $\alpha, \beta>0$, $\alpha_1,
\alpha_2>0$, $m>0$, $\gamma_0>0$, $r_2>0$, $\delta>0$, $s\in
(0,1)$, $d\in (0,1)$, $\varepsilon\in (0,1)$.

System (4) is a system of equations with time delay. The
qualitative analysis is done using the methods from [3].

The equilibrium point of system (4) has the coordinates $Y_0, r_0,
K_0, M_0$, where:
$$Y_0\!=\!\ds\f{g}{d}, r_0\!=\![\ds\f{s(1\!-\!d)}{a}Y_0^{1\!-\!\al_1}]^{-\f{1}{\al_2}},
K_0\!=\!\ds\f{s(1\!-\!d)Y_0}{\delta},
M_0\!=mY_0\!+\!\ds\f{\gamma_0}{r_0\!-\!r_2}.\eqno(5)$$

Using the translation
$$Y=x_1+Y_0, r=x_2+r_0, K=x_3+K_0, M=x_4+M_0$$
in (4) and considering the Taylor expansion of the right members
from (4) until the third order, we have:
$$\dot x(t)=Ax(t)+Bx(t-\tau)+F(x(t), x(t-\tau))\eqno(6)$$
where
$$A\!\!=\!\!\left(\!\!\begin{array}{cccc}
\vspace{0.2cm}
\al a_1 & a\alpha\rho_{01} & 0 & 0\\
\vspace{0.2cm}
\beta m & \beta\gamma_0\gamma_1 & 0 & -\beta\\
\vspace{0.2cm} 0 & a\rho_{01} & -\delta &0\\
a_4 & 0 & 0 & 0
\end{array}\!\!\right)\!\!, B\!\!=\!\!\left(\!\!\begin{array}{cccc}
\vspace{0.2cm}
\al b_1 & 0 & 0 & 0\\
\vspace{0.2cm}
0 & 0 & 0 & 0\\
\vspace{0.2cm}
a\rho_{10} & 0 & 0 & 0\\
\vspace{0.2cm} b_4 & 0 & 0 & 0\end{array}\!\!\right)\eqno(7)$$
where $x(t)=(x_1(t), x_2(t), x_3(t), x_4(t))^T$, and
$$\begin{array}{l}
a_1=(s-1)d(1-\varepsilon)-s+a\rho_{10}, b_1=(s-1)d\varepsilon,
a_4=-d(1-\varepsilon), b_4=-d\varepsilon\\
\rho_{01}=-\al_2Y_0^{\al_1}r_0^{-(\al_2+1)},
\rho_{10}=\al_1Y_0^{\al_1-1}r_0^{-\al_2}\end{array}$$ and

$$F(x(t), x(t-\tau))=(F^1(x(t), x(t-\tau)), F^2(x(t), x(t-\tau)),
F^3(x(t), x(t-\tau)),\eqno(8)$$ $F^4(x(t), x(t-\tau)))^T$,

$ F^1(x(t), x(t-\tau))= \ds\f{\alpha a}{2}\rho_{20}x_1^2(t)+\alpha
a\rho_{11}x_1(t)x_2(t)+\ds\f{\alpha a}{2}\rho_{02}x_2^2(t)+$

           $\ds\f{\alpha a}{6}(\rho_{30}x_1^3(t)+3\rho_{21}x_1(t)^2x_2(t)+3\rho_{12}x_1(t)x_2(t)^2+\rho_{03}x_2(t)^3)$

$$F^2(x(t), x(t-\tau))= \ds\f{\gamma_0}{2}\gamma_2x_2(t)^2
+\ds\f{\gamma_0}{6}\gamma_3x_2(t)^3$$

 $F^3(x(t), x(t-\tau))=
\ds\f{a}{2}\rho_{20}x_1(t-\tau)^2+a\rho_{11}x_1(t-\tau)x_2(t)+\ds\f{a}{2}\rho_{02}x_2(t)^2+$

$+\ds\f{a}{6}(\rho_{30}x_1^3(t-\tau)+3\rho_{21}x_1(t-\tau)^2x_2(t)+3\rho_{12}x_1(t-\tau)x_2(t)^2+\rho_{03}x_2(t)^3)$

$F^4(x(t), x(t-\tau))=0$

and

$$ \rho_{20}=\al_1(\al_1-1)Y_0^{\al_1-2}r_0^{-\al_2}, \quad
\rho_{11}=-\al_1\al_2Y_0^{\al_1-1}r_0^{-(\al_2+1)}$$
$$ \rho_{02}=\al_2(\al_2+1)Y_0^{\al_1}r_0^{-(\al_2+2)}, \quad
\rho_{30}=\al_1(\al_1-1)(\al_1-2)Y_0^{\al_1-3}r_0^{-\al_2}$$
$$ \rho_{21}=-\al_1\al_2(\al_1-1)Y_0^{\al_1-2}r_0^{-(\al_2+1)},
\quad \rho_{12}=\al_1\al_2(\al_2+1)Y_0^{\al_1-1}r_0^{-(\al_2+2)}$$
$$ \rho_{03}=-\al_2(\al_2+1)(\al_2+2)Y_0^{\al_1}r_0^{-(\al_2+3)},$$
$$\gamma_1=-\ds\f{1}{(r_0-r_2)^2},
\gamma_2=\ds\f{2}{(r_0-r_2)^3}, \gamma_3=-\ds\f{6}{(r_0-r_2)^4}.$$

    The characteristic equation of linear part from (6) is:
$$det(\lambda I-A-Be^{-\lambda\tau})=(\lambda+\delta)\Delta(\lambda,
\tau)=0\eqno(9)$$ where
$$\Delta(\lambda,
\tau)=P(\lambda)+e^{-\lam\tau}Q(\lam)\eqno(10)$$

$$P(\lam)=\lambda^3+p_2\lambda^2+p_1\lam+p_0, \quad Q(\lam)=q_2\lambda^2+q_1\lam+q_0$$ and
$$\begin{array}{l}
p_2=-(\al a_1+\beta\gamma_0\gamma_1),\quad
p_1=\al\beta(\gamma_0a_1\gamma_1+ma\rho_{01}), \quad p_0=\al a\beta a_4\rho_{01}\\
q_2=-\al b_1, \quad q_1=\al\beta\gamma_0 b_1\gamma_1, \quad
q_0=\al\beta b_4\rho_{01}.\end{array}$$

To investigate the local stability of the equilibrium point, we
begin by considering, as usual, the case without delay ($\tau=0$
). In this case the characteristic polynomial is:
$$(\lam+\delta)(P(\lambda)+Q(\lam))=0$$
hence, according to the Hurwitz criterion, the equilibrium point
is stable if and only if:
$$p_2+q_2>0, \quad (p_1+q_1)(p_2+q_2)>p_0+q_0.$$

When  $\tau>0$, standard results on stability of systems of delay
differential equations postulate that a equilibrium point is
asymptotically stable if an only if all roots of equation (10)
have a negative real part. It is well known that equation (10) is
a transcendental equation which has an infinite number of complex
roots and the some possible roots with positive real part are
finite in number.

We want to obtain the values $\tau_0$  such that the equilibrium
point (5) changes from local asymptotic stability to instability
or vice versa. We need the imaginary solutions of equation
$\Delta(\lam, \tau)=0$. Let $\lam=\pm i\om$ be these solutions and
without loss of generality we assume $\om>0$. We suppose that
$P(i\om)+Q(i\om)\neq0$, for all $\om\in\rr$. The previous
conditions are equivalent to $(p_1+q_1)(p_2+q_2)\neq p_0+q_0$.

A necessary condition to have $\om$ as a solution of
$\Delta(i\om,\tau)=0$ is that $\om$ must be a root of the
following equation:
$$f(\om)=\om^6+a_F\om^4+b_F\om^2+c_F=0\eqno(11)$$
where $a_F=p_2^2-q_2^2-2p_1$, $b_F=2q_0q_2-2p_0p_2-q_1^2+p_1^2$,
$c_F=p_0^2-q_0^2$.

Let $k=-\ds\f{a_F}{3}$ and
$f_D=\ds\f{1}{4}[f(k)]^4+\ds\f{1}{9}[f^{'}(k)]^3$ .

Using the results from [2],  it results:

{\bf Proposition 1.} {\it 1. Let $\varepsilon>\ds\f{1}{2}$. Then
the following cases can be discerned:

(i) If $a_F\geq0$ or $b_F\leq 0$, then equation (11) has only one
real positive root;

(ii) If $a_F<0$, $b_F>0$  and $f_D\leq0$, then equation (11) has
only one real positive root.


 2. If $\varepsilon<\ds\f{1}{2}$,
$a_F<0$, $b_F<0$ and $f_D\leq0$ then equation (11) has two real
positive roots which are distinct if $f_D\neq0$.

3. If $\varepsilon=\ds\f{1}{2}$ and:

(i) $a_F<0$, $b_F=0$ then equation (11) has only one real positive
root;

(ii)$a_F<0$, $b_F>0,$ then equation (11) has two real positive
roots;

(iii)$b_F<0,$ then equation (11) has only one real positive root.}

Also, we have:

{\bf Theorem 1.} {\it  If we suppose that the equilibrium point
$(Y_0, r_0, K_0, M_0)$ is locally asymptotically stable without
time delay, then in conditions of Proposition 1 there exists only
one stability switch.}

{\bf Theorem 2.} {\it If $\tau_0$ is a stability switch and
$f_D\neq0$, then a Hopf bifurcation occurs at $\tau_0$, where
$$\tau_0=\ds\f{1}{\om_0}arctg\left (\ds\f{\om_0(\om_0^4q_2-
\om^2(q_0-q_1p_2+q_2p_1)+q_0p_1-p_0q_1)}{\om_0^4(q_1-q_2p_2)+\om_0^2(q_0p_2-q_1p_1+p_0q_2)
-p_0q_0}\right)$$ and $\om_0$ is a root of (11).}

\section*{\normalsize\bf 3. The normal form for system (4). Cyclical behavior.}
\vspace{0.6cm}

In this section we describe the direction, stability and the
period of the bifurcating periodic solutions of system (4). The
method we use is based on the normal form theory and the center
manifold theorem introduced by Hassard  [4]. Taking into account
the previous section, if $\tau=\tau_0$ then all roots of equation
(9) other than $\pm i\om_0$  have negative real parts, and any
root of equation (9) of the form $\lam(\tau)=\al(\tau)\pm
i\om(\tau)$ satisfies $\al(\tau_0)=0$, $\om(\tau_0)=\om_0$ and
$\ds\f{d\al(\tau_0)}{d\tau}\neq0$. For notational convenience let
$\tau=\tau_0+\mu, \mu\in\rr$. Then $\mu=0$ is the Hopf bifurcation
value for equations (4).

Define the space of continuous real-valued functions as
$C=C([-\tau_0,0],\rr^4).$

In $\tau=\tau_0+\mu, \mu\in\rr$, we regard $\mu$ as the
bifurcation parameter. For $\Phi\in C$ we define a linear
operator:
$$L(\mu)\Phi=A\Phi(0)+B\Phi(-\tau)$$ where A and B are given by
(7) and a nonlinear operator $F(\mu, \Phi)=F(\Phi(0)$,
$\Phi(-\tau))$, where $F(\Phi(0), \Phi(-\tau))$ is given by (8).
By the Riesz representation theorem, there exists a matrix whose
components are bounded variation functions, $\eta(\theta, \mu)$
with $\theta\in[-\tau_0, 0]$ such that:
$$L(\mu)\Phi=\int\limits_{-\tau_0}^0d\eta(\theta,\mu)\phi(\theta),
\quad \theta\in[-\tau_0,0].$$

For $\Phi\in C^1([-\tau_0, 0], \rr^{4})$ we define:
$$\ca(\mu)\Phi(\th)=\left\{\begin{array}{ll} \vspace{0.2cm}
\ds\f{d\Phi(\th)}{d\th}, & \th\in[-\tau_0,0)\\
\int\limits_{-\tau_0}^0d\eta(t,\mu)\phi(t), &
\th=0,\end{array}\right.$$
$$R(\mu)\Phi(\th)=\left\{\begin{array}{ll} \vspace{0.2cm}
0, & \th\in[-\tau_0,0)\\
F(\mu, \Phi), & \th=0.\end{array}\right.$$

We can rewrite (6) in the following vector form:
$$\dot u_t=\ca(\mu)u_t+R(\mu)u_t\eqno(12)$$
where $u=(u_1, u_2, u_3, u_4)^T$, $u_t=u(t+\theta)$ for
$\theta\in[-\tau_0, 0]$.

For $\Psi\in C^1([0,\tau_0], \rr^{*4})$, we define the adjunct
operator $\ca^*$ of $\ca$ by:
$$\ca^*\Psi(s)=\left\{\begin{array}{ll} \vspace{0.2cm}
-\ds\f{d\Psi(s)}{ds}, & s\in(0, \tau_0]\\
\int\limits_{-\tau_0}^0d\eta^T(t,0)\psi(-t), &
s=0.\end{array}\right.$$ We define the following bilinear form:
$$<\Psi(\th),
\Phi(\th)>=\bar\Psi^T(0)\Phi(0)-\int_{-\tau_0}^0\int_{\xi=0}^0\bar\Psi^T(\xi-\theta)d\eta(\theta)\phi(\xi)d\xi,$$
where $\eta(\theta)=\eta(\theta,0)$.

We assume that $\pm i\om_0$ are eigenvalues of $\ca(0)$. Thus,
they are also eigenvalues of $\ca^*$. We can easily obtain:
$$\Phi(\th)=ve^{\lam_1\th},\quad \th\in[-\tau_0, 0]\eqno(13)$$
where $v=(v_1, v_2, v_3, v_4)^T$,
$$v_1=1, v_2=-\ds\f{\beta(a_4+b_4e^{-\lam_1\tau_0}-m\lam_1)}{\lam_1(\lam_1-\beta\gamma_0\gamma_1)}, v_4=\ds\f{a_4+b_4e^{-\lam_1\tau_0}}
{\lam_1},$$
$$v_3=\ds\f{a}{\lam_1+\delta}[\rho_{10}e^{-\lam_1\tau_0}+\rho_{01}\ds\f{\beta(a_4+b_4e^{-\lam_1\tau_0}-m\lam_1)}{\lam_1(\beta\gamma_0\gamma_1-\lam_1)}]$$
is the eigenvector of $\ca(0)$ corresponding to $\lam_1=i\om_0$
and
$$\Psi(s)=we^{\lam_2s},\quad s\in[0,\infty)$$ where
$w=(w_1, w_2, w_3, w_4)$,
$$w_1=\ds\f{\lam_2-\beta\gamma_0\gamma_1}{\alpha a\rho_{01}}\ds\f{1}{\bar\eta}, w_2=\ds\f{1}{\bar\eta},
w_3=0, w_4=-\ds\f{\beta}{\lam_2}\ds\f{1}{\bar\eta}$$
$$\eta=\ds\f{\lam_1-\beta\gamma_0\gamma_1}{a\alpha\rho_{01}}(1+\al b_1\ds\f{\lam_1\tau_0e^{-\lam_1\tau_0}\!-\!e^{-\lam_1\tau_0}+1}{\lam_1^2})+v_2-
\ds\f{\beta}{\lam_1}
(v_4+b_4\ds\f{\lam_1\tau_0e^{-\lam_1\tau_0}\!-\!e^{-\lam_1\tau_0}+1}{\lam_1^2})$$
is the eigenvector of $\ca(0)$ corresponding to $\lam_2=-i\om_0.$

We can verify that: $<\Psi(s), \Phi(s)>=1$, $<\Psi(s),
\bar\Phi(s)>=<\bar\Psi(s), \Phi(s)>=0$, $<\bar\Psi(s),
\bar\Phi(s)>=1.$

Using the approach of Hassard [4], we next compute the coordinates
to describe the center manifold $\Omega_0$ at $\mu=0$. Let
$u_t=u_t(t+\th), \th\in[-\tau_0,0)$ be the solution of equation
(12) when $\mu=0$ and
$$z(t)=<\Psi, u_t>,
\quad w(t,\th)=u_t(\th)-2Re\{z(t)\Phi(\th)\}.$$

On the center manifold $\Omega_0$, we have:
$$w(t,\th)=w(z(t), \bar z(t), \th)$$ where
$$w(z,\bar z, \th)=w_{20}(\th)\ds\f{z^2}{2}+w_{11}(\th)z\bar
z+w_{02}(\th)\ds\f{\bar z^2}{2}+w_{30}(\th)\ds\f{z^3}{6}+\dots$$
in which $z$ and $\bar z$ are local coordinates for the center
manifold $\Omega_0$ in the direction of $\Psi$ and $\bar\Psi$ and
$w_{02}(\th)=\bar w_{20}(\th)$. Note that $w$ and $u_t$ are real.

For solution $u_t\in \Omega_0$ of equation (12), since $\mu=0$, we
have:
$$\dot z(t)=\lam_1z(t)+g(z, \bar z)\eqno(14)$$ where
$$\begin{array}{ll}
g(z, \bar z)& =\bar\Psi(0)F(w(z(t),\bar z(t), 0)+2Re(z(t)\Phi(0)))=\\
& =g_{20}\ds\f{z(t)^2}{2}+g_{11}z(t)\bar z(t)+g_{02}\ds\f{\bar
z(t)^2}{2}+g_{21}\ds\f{z(t)^2\bar z(t)}{2}+\dots\end{array}$$
where
$$g_{20}=F^1_{20}\bar w_1+F^2_{20}\bar w_2,
g_{11}=F^1_{11}\bar w_1+F^2_{11}\bar w_2, g_{02}=F^1_{02}\bar
w_1+F^2_{02}\bar w_2,\eqno(15)$$ with
$$\begin{array}{l}F_{20}^1=\alpha a(\rho_{20}+2\rho_{11}v_2+\rho_{02}v_2^2), F_{11}^1=\alpha a(\rho_{20}+\rho_{11}(\bar v_2+v_2)+\rho_{02}v_2\bar v_2), \\
F_{20}^2=\gamma_0\gamma_2v_2^2, F_{11}^2=\gamma_0\gamma_2v_2\bar
v_2, F_{02}^1= \bar F_{20}^1, F_{02}^2=\bar F_{20}^2,
\end{array}$$ and
$$g_{21}=F^1_{21}\bar w_1+F^2_{21}\bar w_2\eqno(16)$$ where
$$\begin{array}{l}
F_{21}^1= \alpha a\rho_{20}(2w_{11}^1(0)+w_{20}^1(0))+\alpha
a\rho_{11}(2w_{11}^2(0)+w_{20}^2(0)+2w_{11}^1(0)\bar
v_2 +\\
w_{20}^1(0)v_2)+ \alpha a\rho_{02}(2w_{11}^2(0)v_2+w_{20}^2(0)\bar
v_2)+\alpha
a(\rho_{30}+2\rho_{21}v_2+2\rho_{12}v_2^2+\\
+\rho_{03}v_2^2\bar
v_2 +\rho_{12}v_2^2+\rho_{21}\bar v_2)\end{array}$$
$$F_{21}^2=\gamma_0\gamma_2(2w_{11}^2(0)v_2+w_{20}^2(0)\bar
v_2)+\gamma_0\gamma_3v_2^2\bar v_2.$$

The vectors $w_{20}(\theta)$, $w_{11}(\theta)$ with
$\theta\in[-\tau,0]$ are given by:
$$\begin{array}{l}
w_{20}(\theta)=-\ds\f{g_{20}}{\lam_1}ve^{\lam_1\theta}-\ds\f{\bar
g_{02}}{3\lam_1}\bar v_e^{\lam_2\theta}+E_1e^{2\lam_1\theta}\\
w_{11}(\theta)=\ds\f{g_{11}}{\lam_1}ve^{\lam_1\theta}-\ds\f{\bar
g_{11}}{\lam_1}\bar ve^{\lam_2\theta}+E_2\end{array}\eqno(17)$$
where
$$
E_1=-(A+e^{-\lam_1\tau_0}B-2\lam_1I)^{-1}F_{20},\quad
E_2=-(A+B)^{-1}F_{11},$$ where $F_{20}=(F_{20}^1, F_{20}^2,
F_{20}^3, 0)^T$, $F_{11}=(F_{11}^1, F_{11}^2, F_{11}^3, 0)^T$.

Based on the above analysis and calculation, we can see that each
$g_{ij}$ in (15), (16) are determined by the parameters and delay
from system (4). Thus, we can explicitly compute the following
quantities:
$$\begin{array}{l}
C_1(0)=\ds\f{i}{2\om_0}(g_{20}g_{11}-2|g_{11}|^2-\ds\f{1}{3}|g_{02}|^2)+\ds\f{g_{21}}{2}\\
\mu_2=-\ds\f{Re(C_1(0))}{Re\lam'(\tau_0)},
T_2=-\ds\f{Im(C_1(0))+\mu_2Im\lam'(\tau_0)}{\om_0},
\beta_2=2Re(C_1(0)).\end{array}\eqno(18)$$ and
$$\lam'=\ds\f{\lam(q_2\lam^2+q_1\lam+q_0)e^{-\lam \tau}}
{3\lam^2+2p_2\lam+p_1+e^{-\lam\tau}(2q_2\lam+q_1-\tau(q_2\lam^2+q_1\lam+q_0))}$$
In summary, this leads to the following result:

\vspace{2mm} {\bf Theorem 3.} {\it In formulas (18), $\mu_2$
determines the direction of the Hopf bifurcation: if $\mu_2>0
(<0)$, then the Hopf bifurcation is supercritical (subcritical)
and the bifurcating periodic solutions exit for $\tau>\tau_0
(<\tau_0)$; $\beta_2$ determines the stability of the bifurcating
periodic solutions: the solutions are orbitally stable (unstable)
if $\beta_2<0 (>0)$; and $T_2$ determines the period of the
bifurcating periodic solutions: the period increases (decreases)
if $T_2>0 (<0)$.}

\section*{\normalsize\bf 4. Numerical example.}
\vspace{0.6cm}

\hspace{0.6cm} In this section we find the waveform plots through
the formula:
$$X(t+\th)\!=\! z(t)\Phi(\th)\!+\!\bar
z(t)\bar\Phi(\th)\!+\!\ds\f{1}{2}w_{20}(\th)z^2(t)+w_{11}(\th)z(t)\bar
z(t)\!+\!\ds\f{1}{2}w_{02}(\th)\bar z(t)^2+X_0,$$ where $z(t)$ is
the solution of (14), $\Phi(\th)$ is given by (13), $w_{20}(\th),
w_{11}(\th), w_{02}(\th)=\bar w_{20}(\th)$ are given by (17) and
$X_0=(Y_0, r_0, K_0, M_0)^T$ is the equilibrium state.

For the numerical simulations we use Maple 9.5. We consider system
(4) with $a=0.38$, $\alpha=0.96$, $\beta=1$, $\alpha_1=0.5$,
$\alpha_2=0.83$, $\gamma_0=1$, $d=0.1$, $s=0.3$, $r_2=0.003$,
$\delta=0.2$, $m=0.005$, $g=50$. The equilibrium point is:
$Y_0\!=\! 500$, $r_0\!=\! 0.03572181612$, $K_0\!=\! 675$,
$M_0\!=\! 33.06065092$. In what follows we consider two different
shares $\epsilon$ of delay tax revenues: $\varepsilon=0.3$ and
$\varepsilon=0.8$.

For $\varepsilon=0.3$ we obtain: $\mu_2\!=\!
1.654628706\cdot10^{-8}$, $\beta_2\!=\! 2.224294680\cdot10^{-9}$,
$T_2\!=\! 2.092652051\cdot10^{-9}$, $\omega_0\!=\! 0.6685954740$,
$\tau_0\!=\! 4.965007916$. Then the Hopf bifurcation is
supercri\-ti\-cal, the solutions are orbitally unstable and the
period of the solution is increasing. The wave plots are given in
the following figures:

\begin{center}
{\small \begin{tabular}{c|c|c} \hline
Fig.1.Waveplot $(t,Y(t))$&Fig.2.Waveplot $(t,r(t))$&Fig.3.Waveplot $(t,K(t))$\\&&\\
\cline{1-3} \epsfxsize=4cm

\epsfysize=4cm

\epsffile{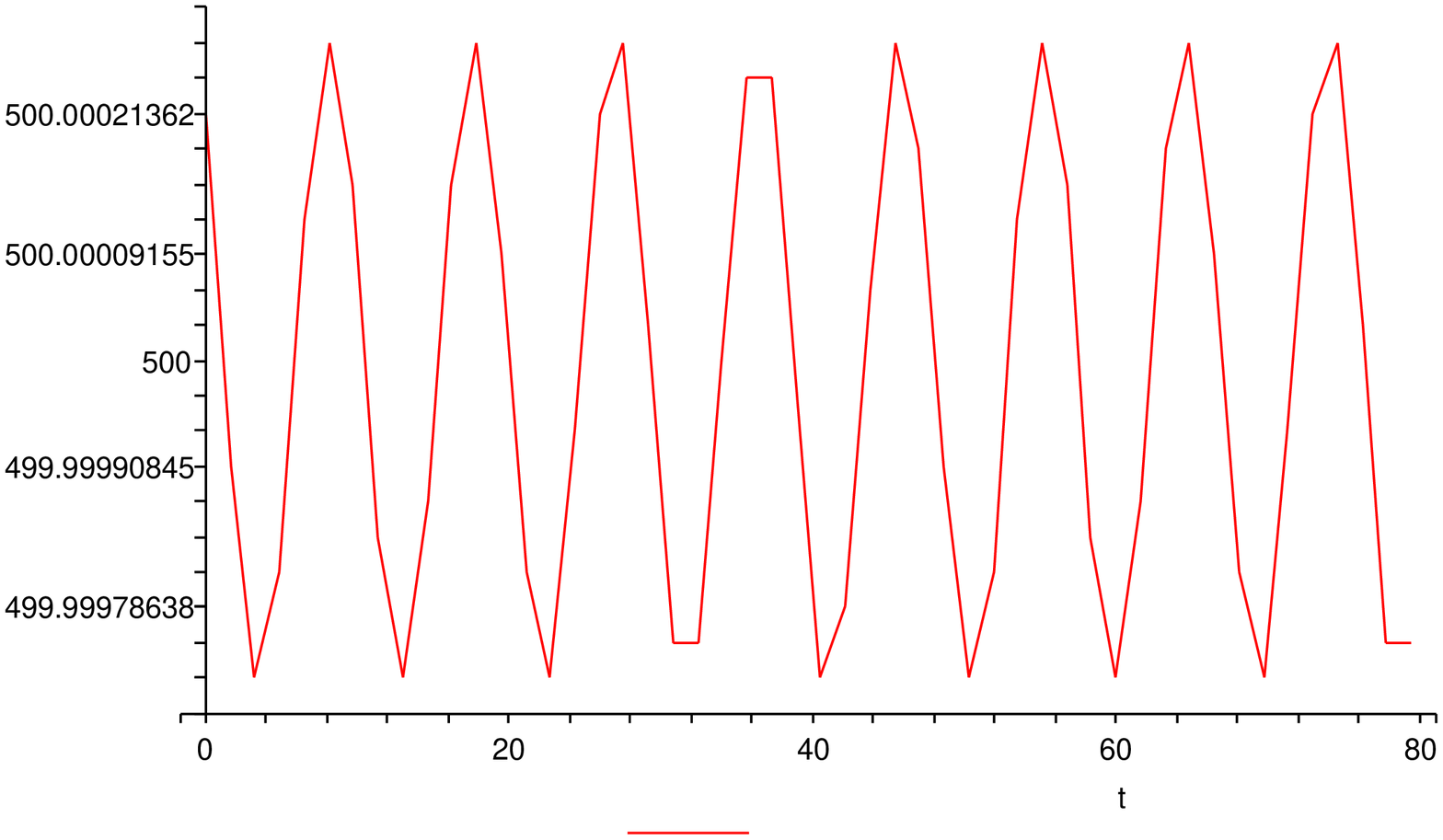} &

\epsfxsize=4cm

\epsfysize=4cm

\epsffile{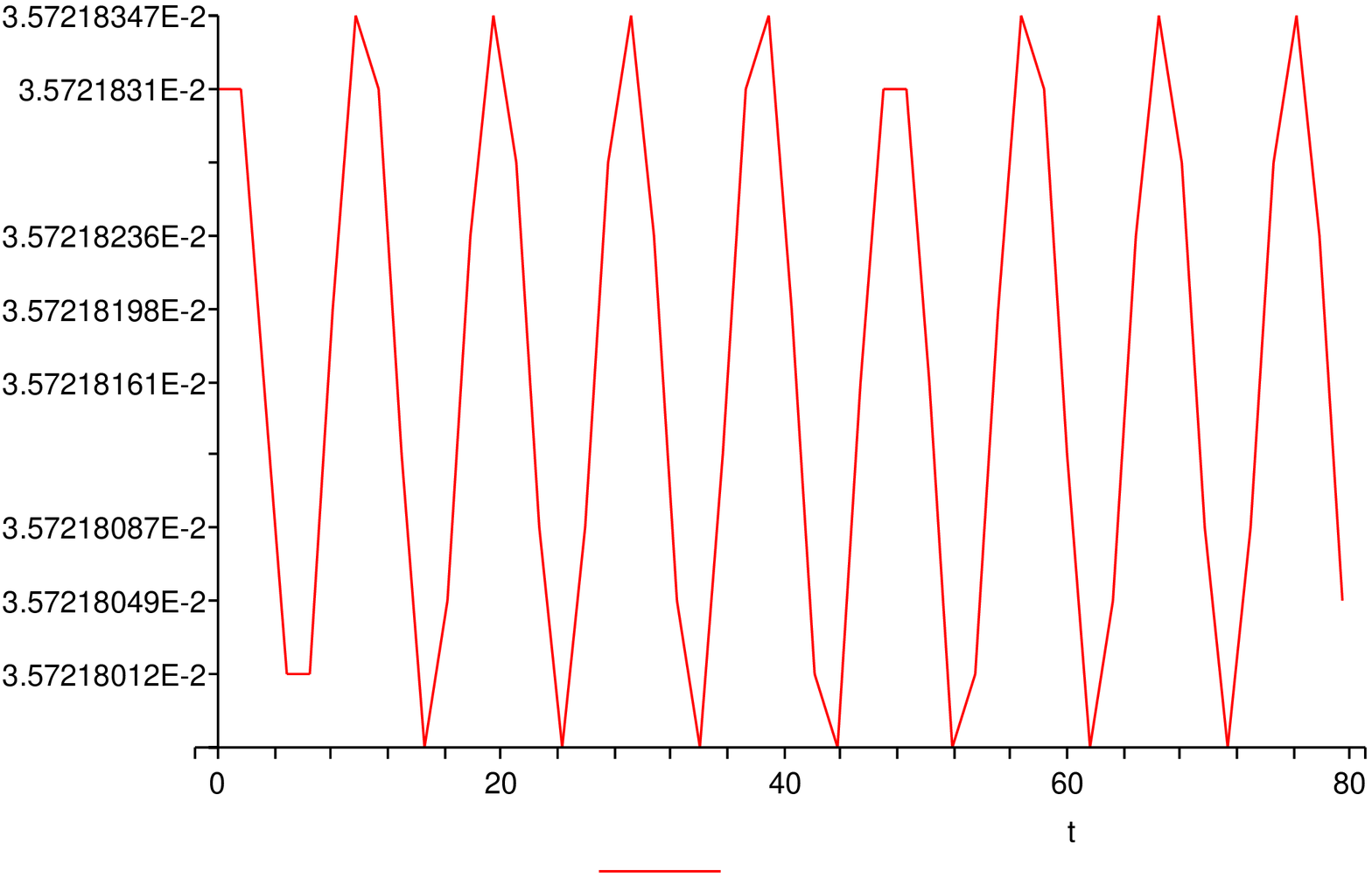} &

\epsfxsize=4cm

\epsfysize=4cm

\epsffile{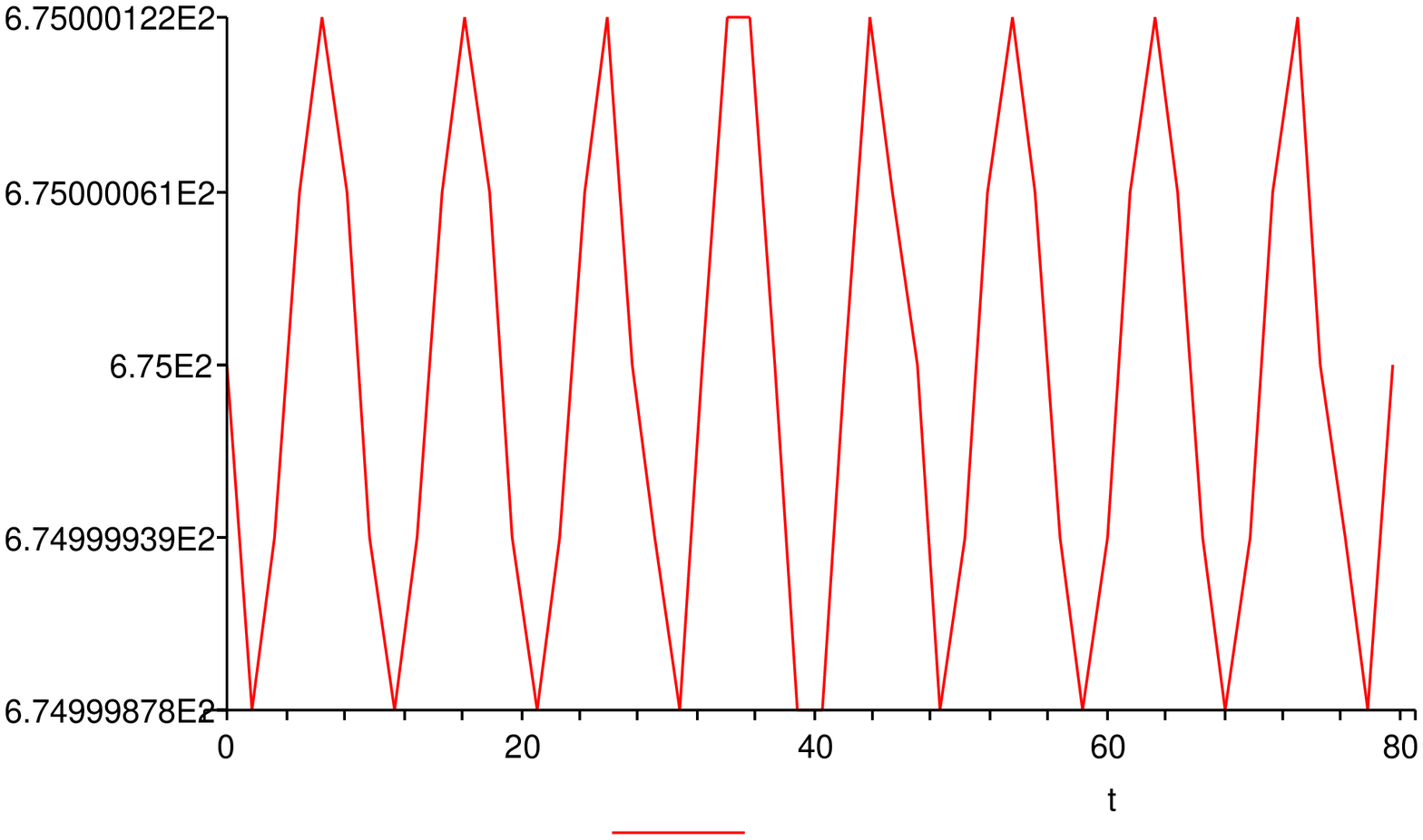}
\\
 \hline
\end{tabular}}
\end{center}

\medskip

\begin{center}
{\small \begin{tabular}{c|c|c} \hline
Fig.4.Waveplot $(t,M(t))$&Fig.5.Waveplot $(t,I(t))$&Fig.6.Waveplot $(t,L(t))$\\&&\\
\cline{1-3} \epsfxsize=4cm

\epsfysize=4cm

\epsffile{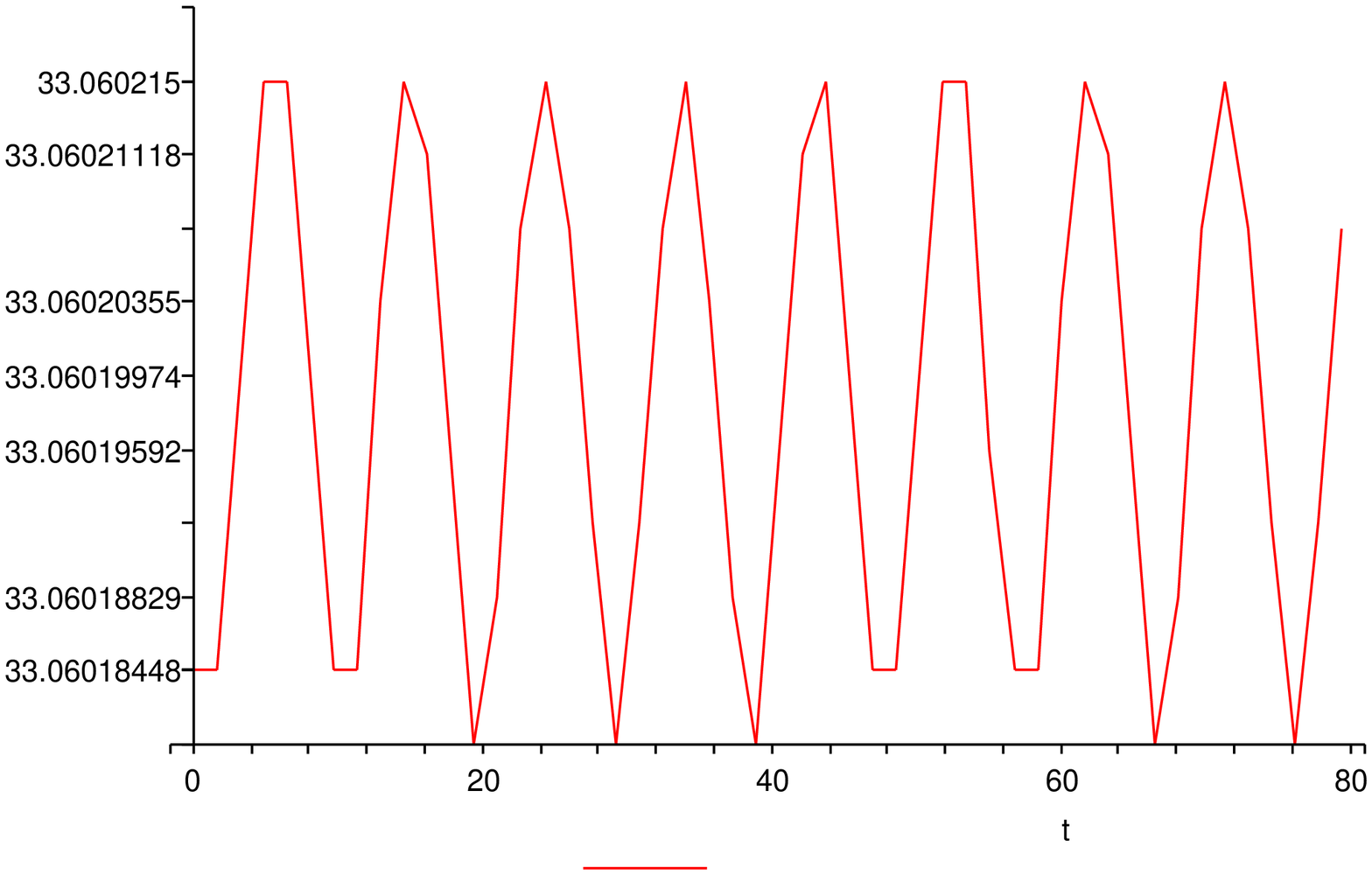} &

\epsfxsize=4cm

\epsfysize=4cm

\epsffile{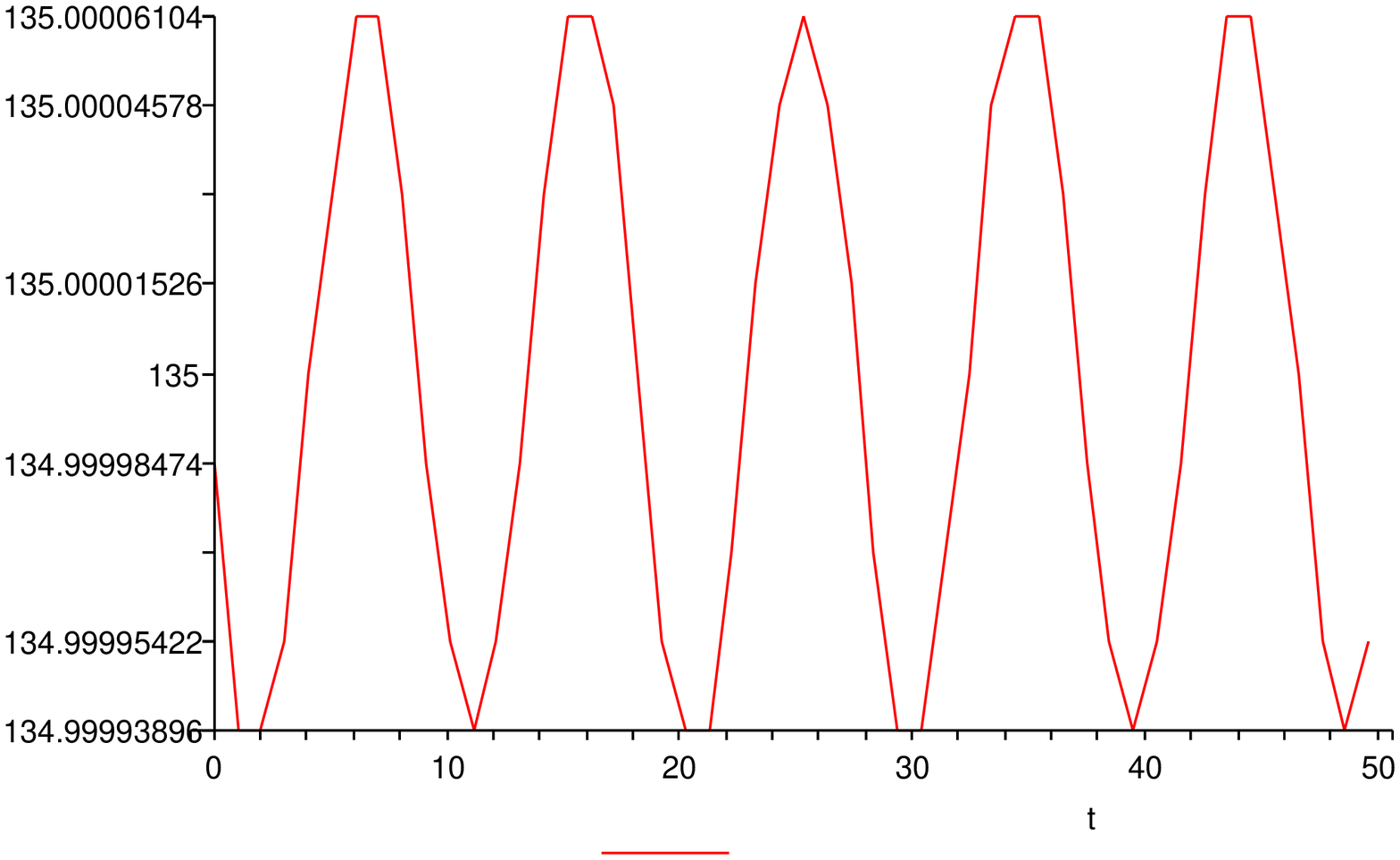} &

\epsfxsize=4cm

\epsfysize=4cm

\epsffile{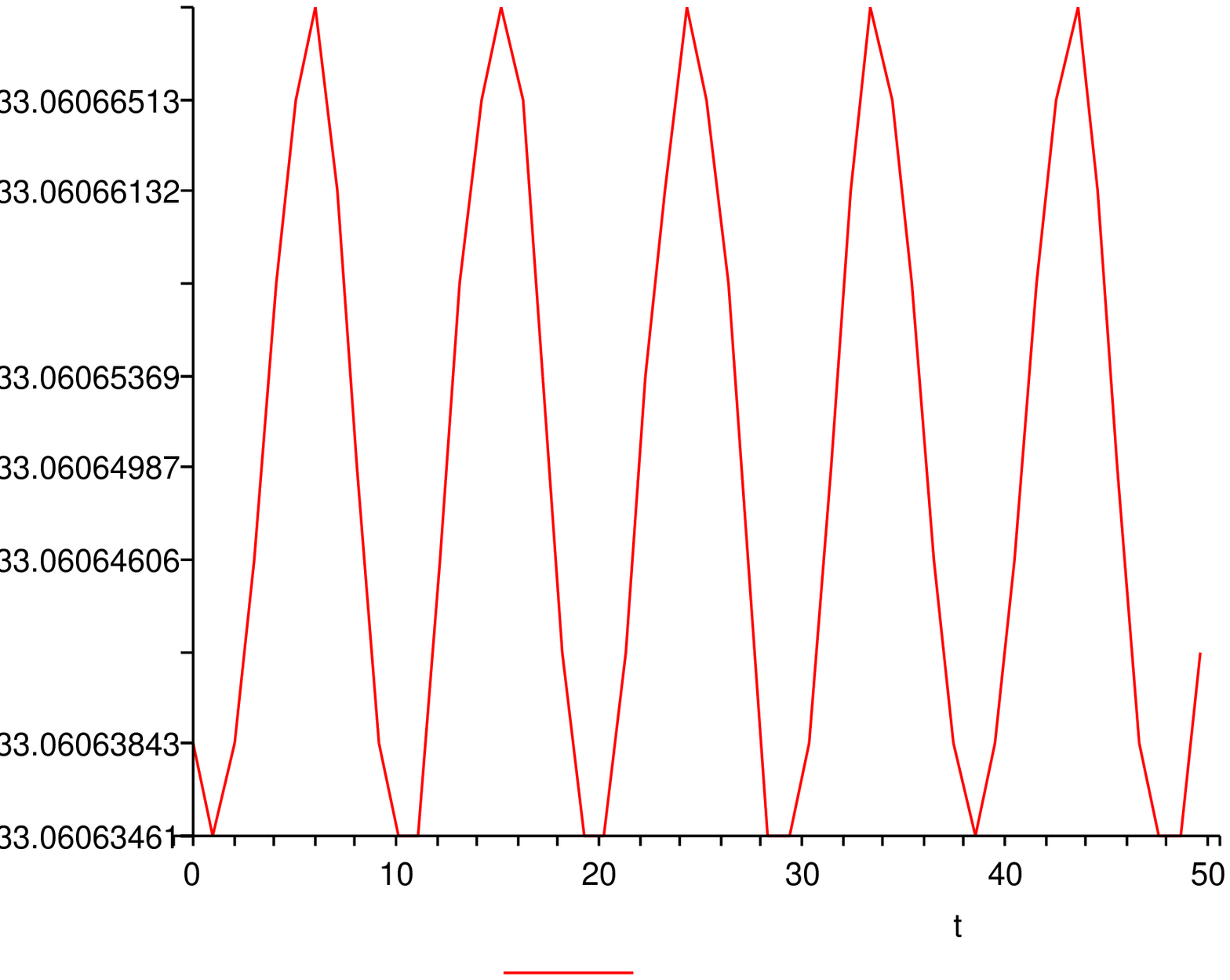}
\\
 \hline
\end{tabular}}
\end{center}

\medskip

\medskip
\medskip

For $\varepsilon=0.8$ we obtain: $\mu_2\!=\!
-9.160756314\cdot10^{-8}$, $\beta_2\!=\!
-1.968119398\cdot10^{-8}$, $T_2\!=\! -1.608068638\cdot10^{-8}$,
$\omega_0\!=\! 0.8553440397$, $\tau_0\!=\! 3.918246696$. Then the
Hopf bifurcation is subcritical, the solutions are orbitally
stable and the period of the solution is decreasing. The wave
plots are given in the following figures:

\begin{center}
{\small \begin{tabular}{c|c|c} \hline
Fig.7.Waveplot $(t,Y(t))$&Fig.8.Waveplot $(t,r(t))$&Fig.9.Waveplot $(t,K(t))$\\&&\\
\cline{1-3} \epsfxsize=4cm

\epsfysize=4cm

\epsffile{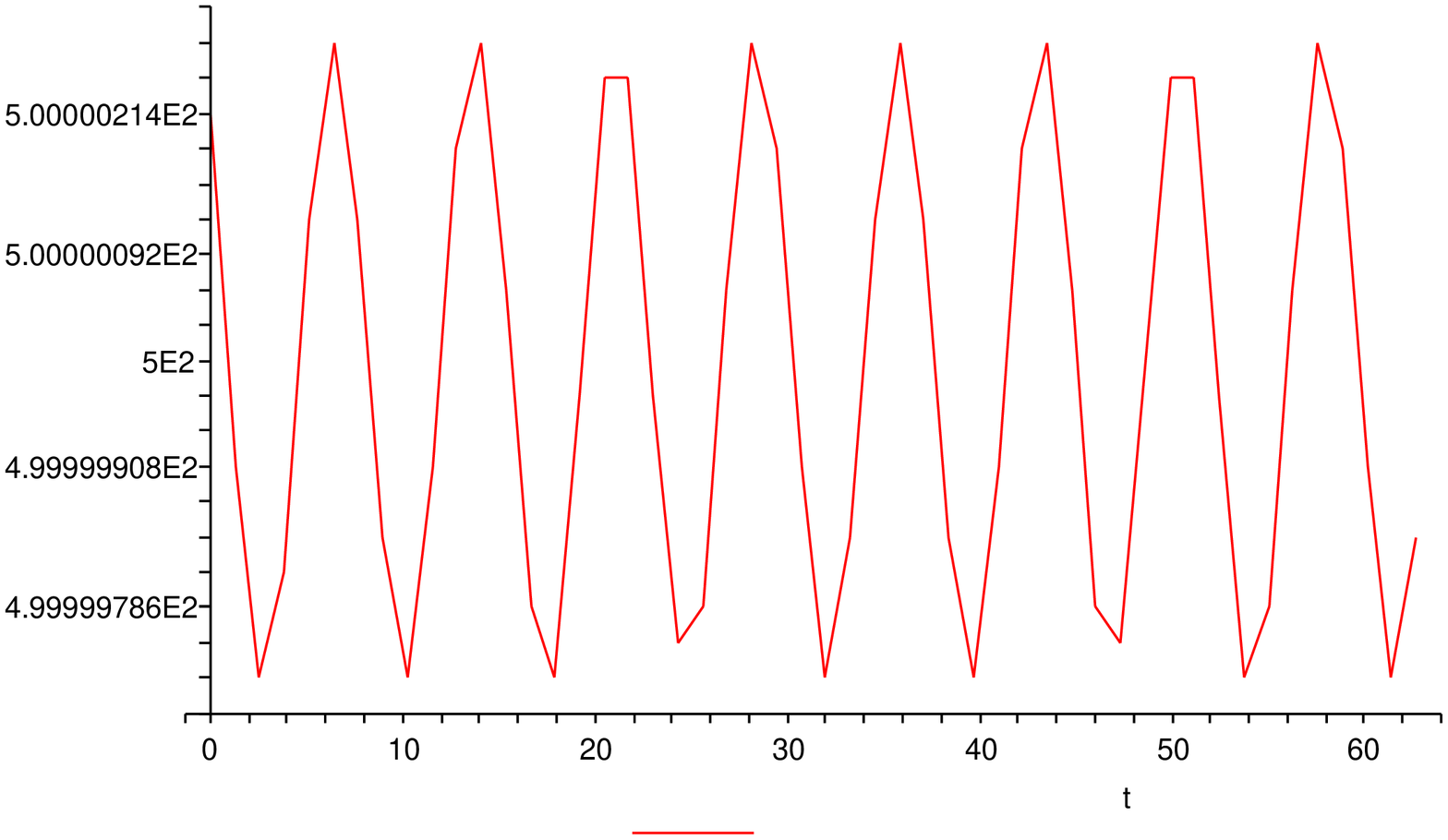} &

\epsfxsize=4cm

\epsfysize=4cm

\epsffile{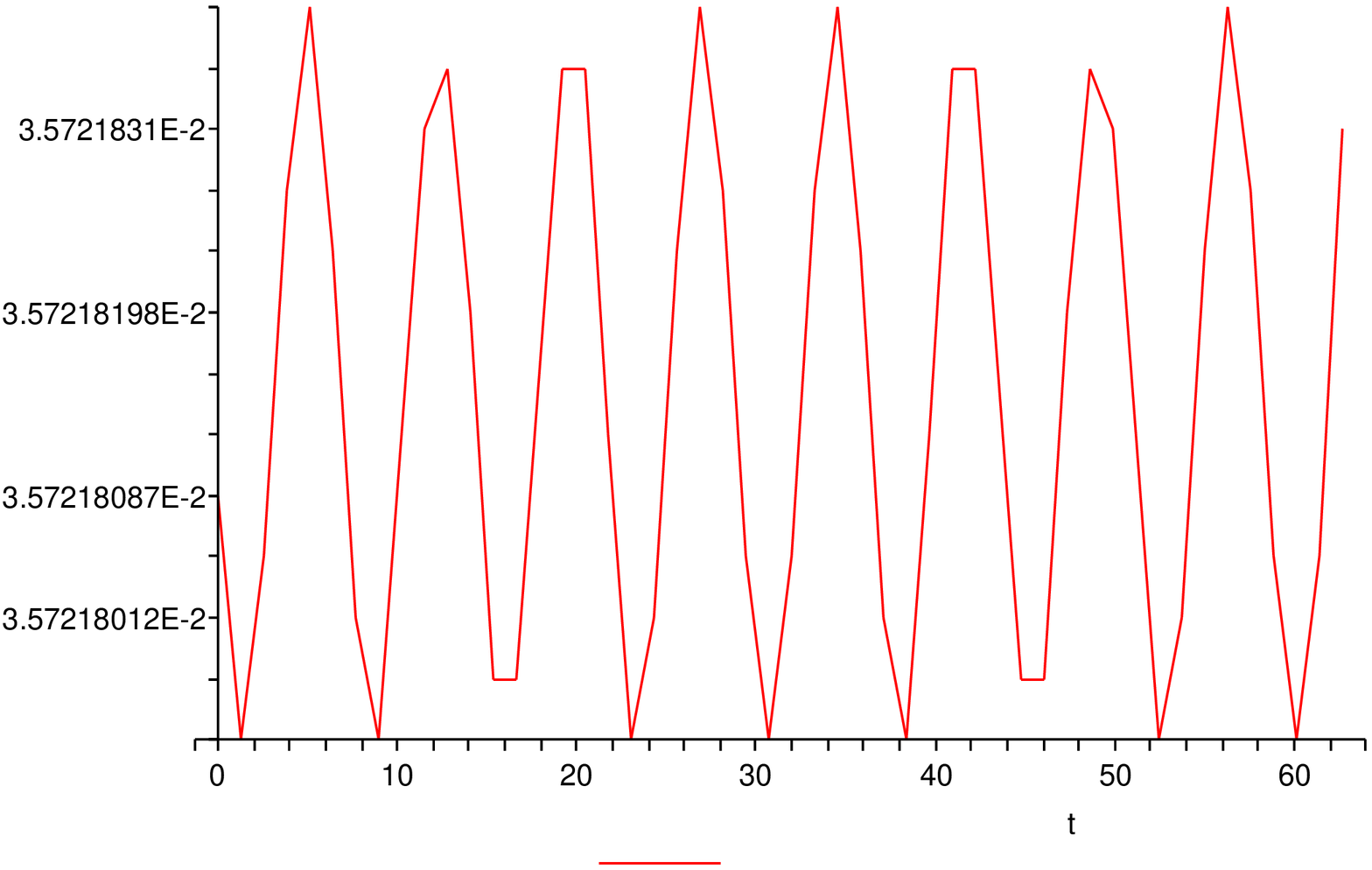} &

\epsfxsize=4cm

\epsfysize=4cm

\epsffile{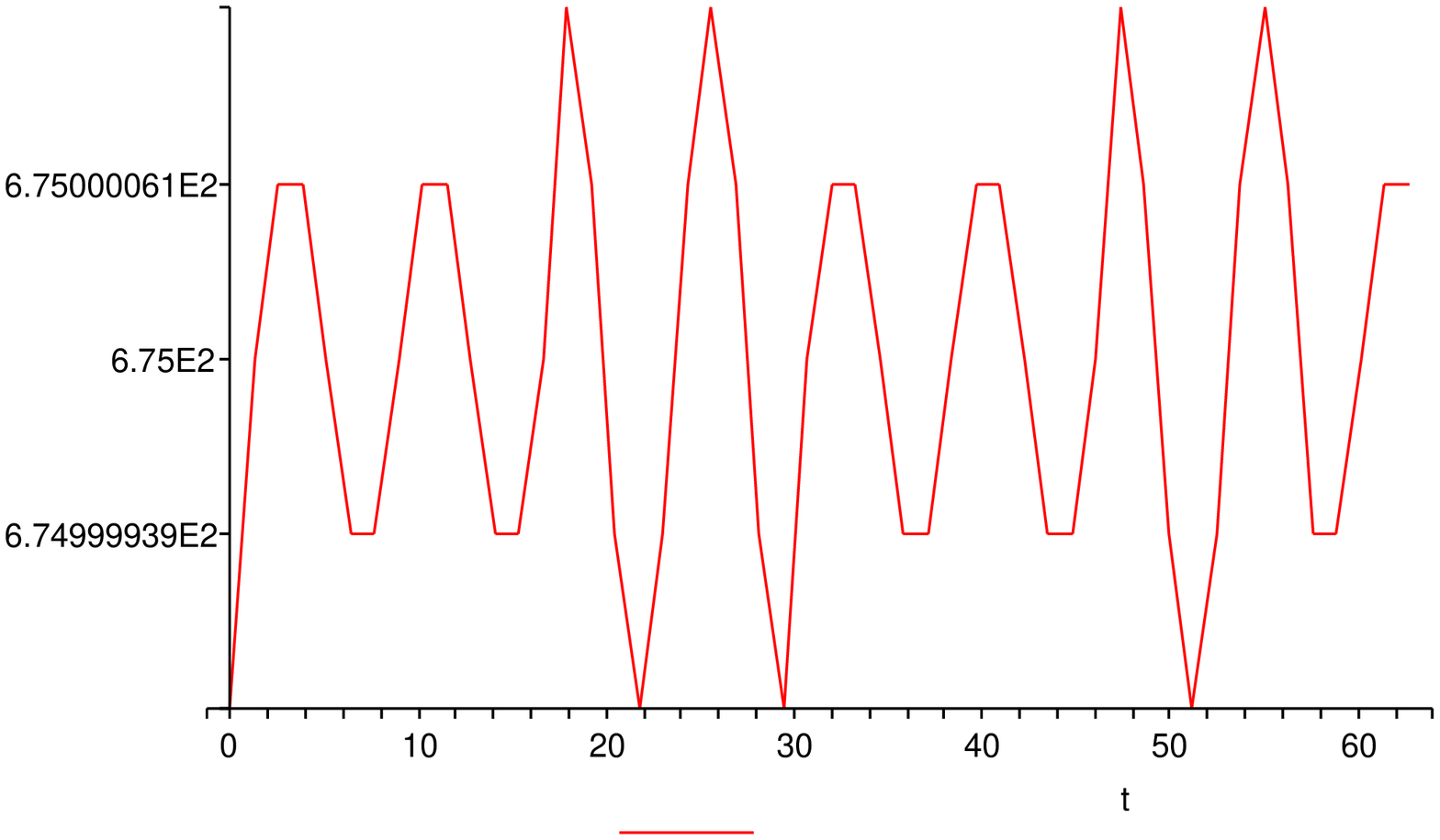}
\\
 \hline
\end{tabular}}
\end{center}

\medskip

\begin{center}
{\small \begin{tabular}{c|c|c} \hline
Fig.10.Waveplot $(t,M(t))$&Fig.11.Waveplot $(t,I(t))$&Fig.12.Waveplot $(t,L(t))$\\&&\\
\cline{1-3} \epsfxsize=4cm

\epsfysize=4cm

\epsffile{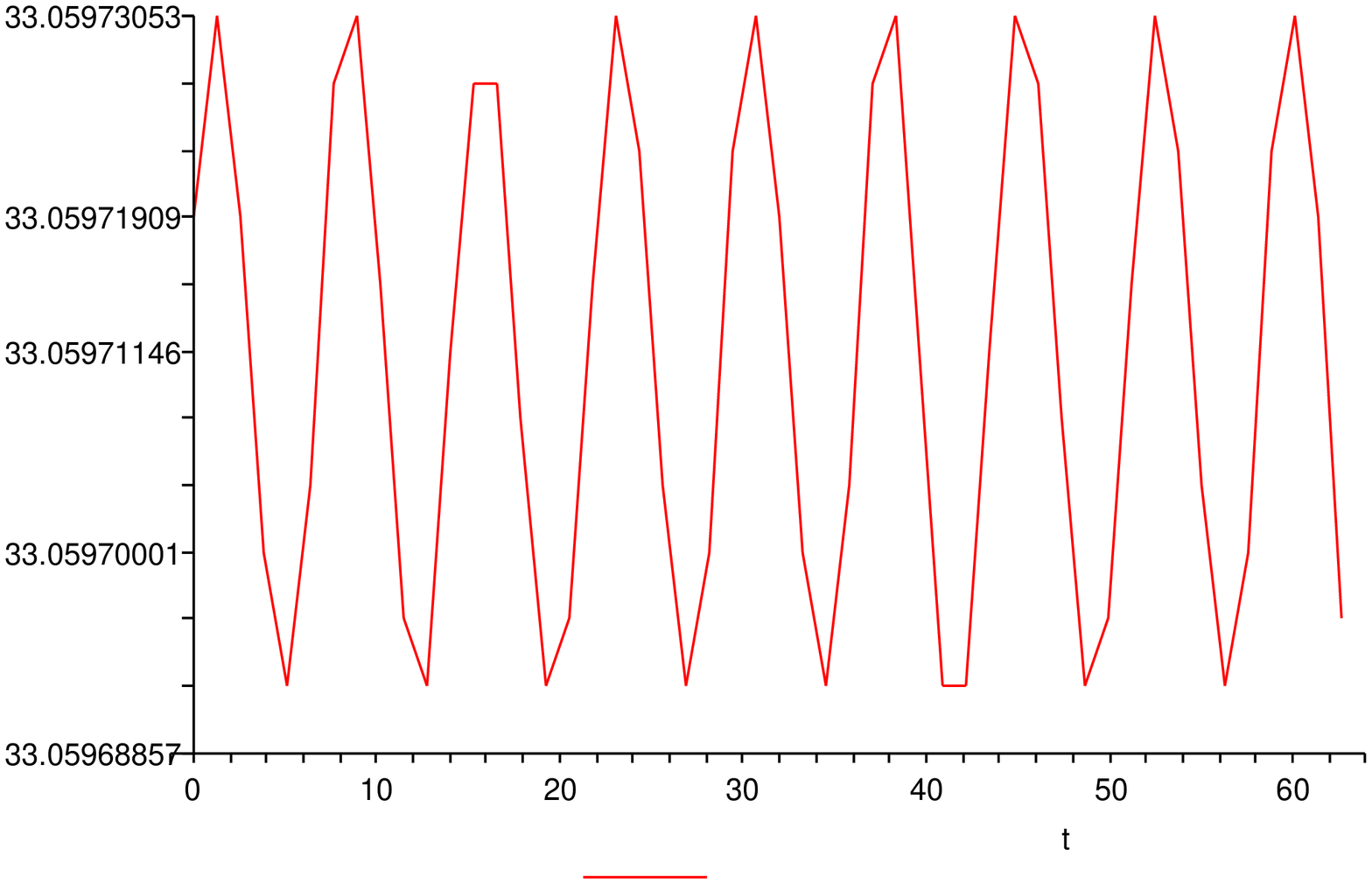} &

\epsfxsize=4cm

\epsfysize=4cm

\epsffile{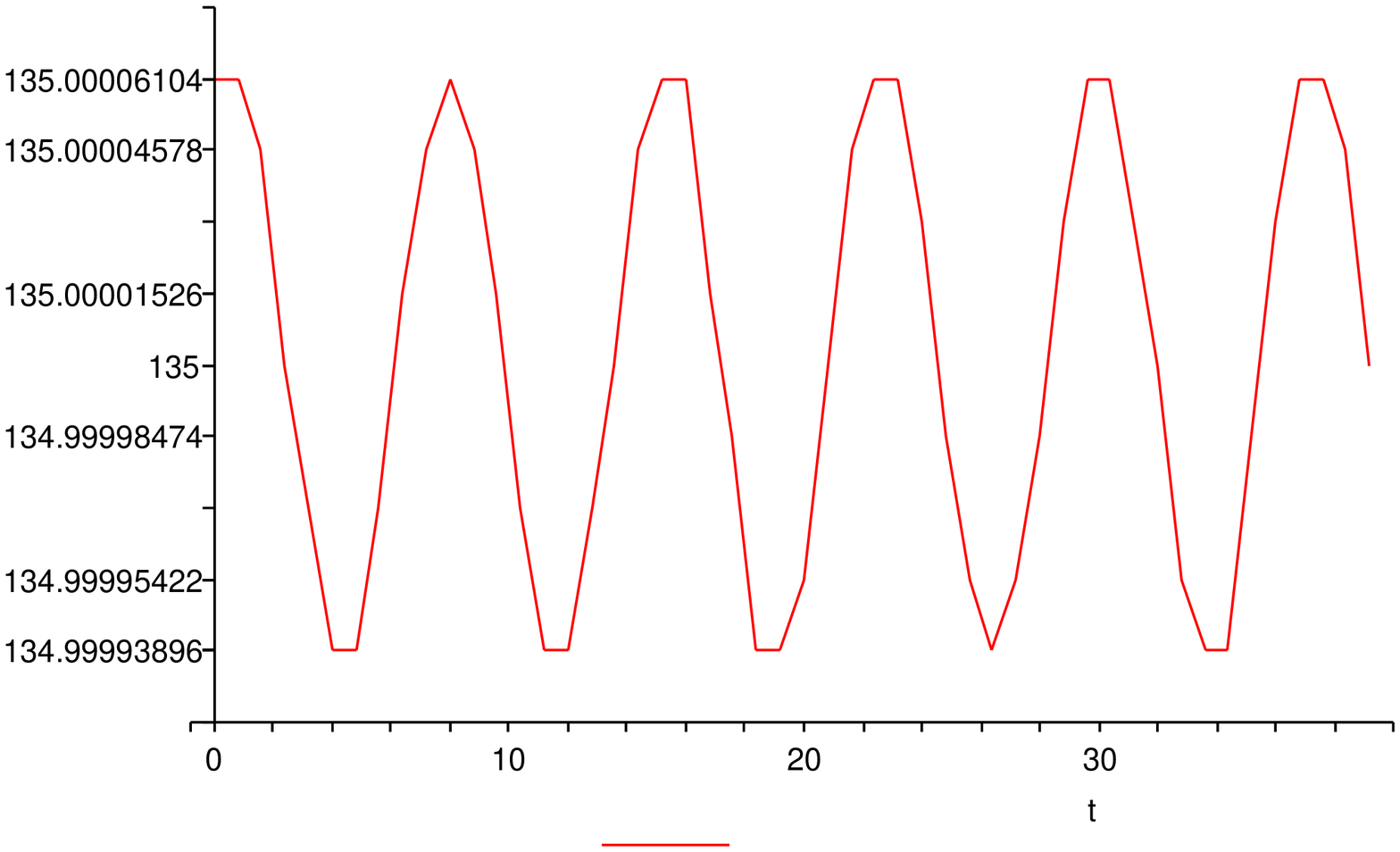} &

\epsfxsize=4cm

\epsfysize=4cm

\epsffile{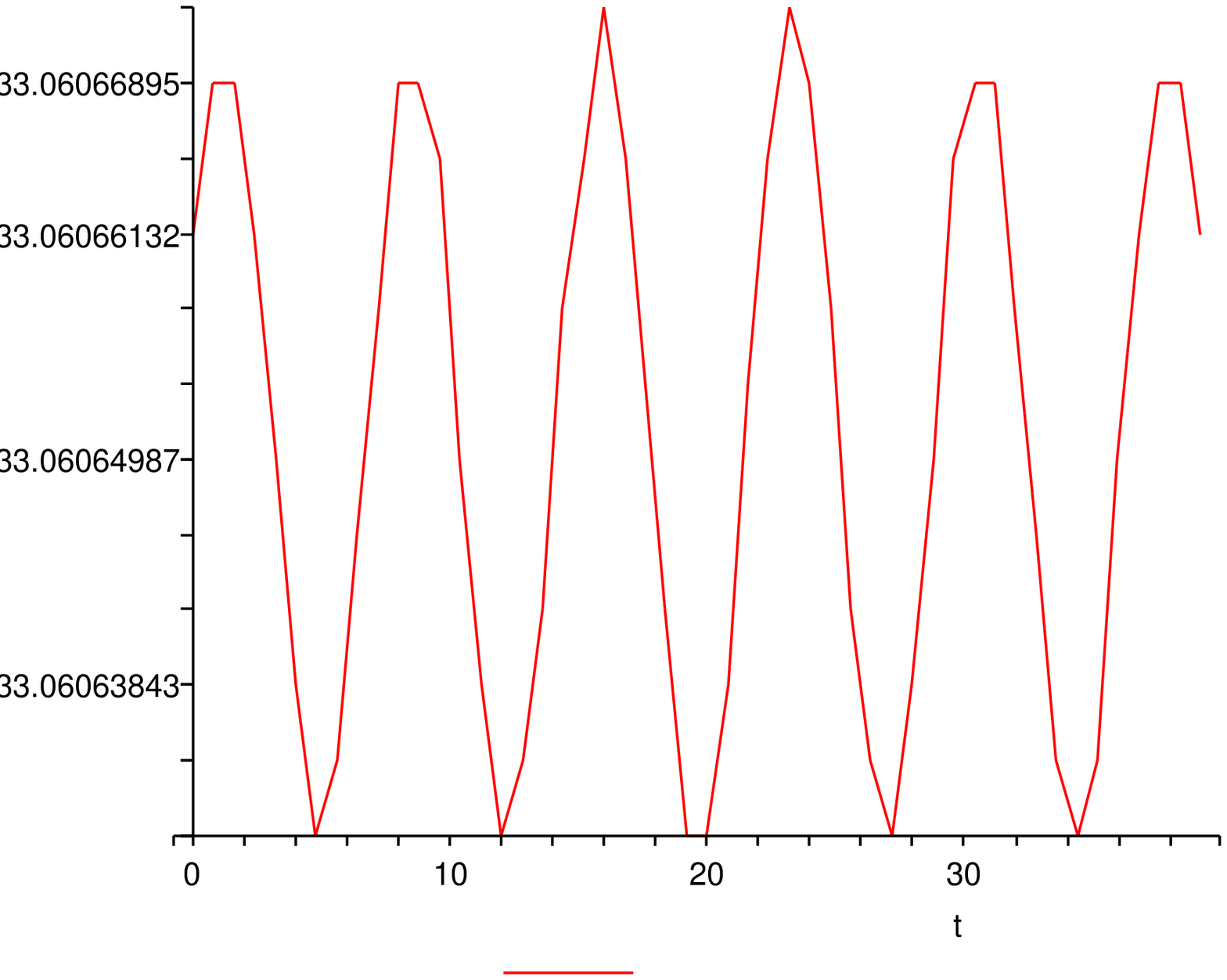}
\\
 \hline
\end{tabular}}
\end{center}

\medskip

\medskip

\section*{\normalsize\bf 5. Conclusions.}
\vspace{0.6cm}

From the analysis of the model with continuous time, it results
that the model accepts a limit cycle. The nature of the limit
cycle is given by the coefficients (18) which include the
parameters of the model. We establish the nature of the limit
cycle. Because the expressions of the parameter and the
coefficients (18) are difficult to analyze directly, the use of
Maple 9.5 was essential. The paper's results confirm that a series
of economical processes, where a variable with time delay
intervenes, have a limit cycle, thus, allowing a prediction
concerning the evolution of the model.

\end{document}